\documentclass[a4paper,12pt]{article}
\usepackage{latexsym,amssymb,amsmath,amsthm}
\usepackage{enumerate}
\usepackage{color}
\usepackage{graphicx}
\usepackage[all,poly]{xy}

\newtheorem{thm}{Theorem}[section]
\newtheorem{prop}[thm]{Proposition}
\newtheorem{lemma}[thm]{Lemma}

\numberwithin{equation}{thm}
\numberwithin{thm}{section}

\newcommand{\ran}{{\rm ran}\,}

\newcommand{\C}{\ensuremath{\mathcal{C}}}

\newcommand{\K}{\ensuremath{\mathcal{K}}}
\def \b{\beta}
\def \a{\alpha}
\def \v{\varphi}
\def \t{\theta}

\begin{document}
\begin{center}

\textbf{\large  A note on isomorphism theorems for semigroups of order-preserving transformations with restricted range}

\end{center}

\centerline {   Phichet Jitjankarn $^\dag$   \ \  and \ \   Thitarie Rungratgasame {\footnote{This research was partially supported by SWU endownment fund year 2014.}}      
}
\renewcommand{\thefootnote}{\arabic{footnote}}

\vskip.5cm
\centerline{$^\dag$Division of Mathematics, School of Science, Walailak University, Nakhon Si Thammarat 80161,
Thailand
% address goes in here
}
\centerline{\texttt{$^\dag$jitjankarn@gmail.com}
% email goes in here
}
\centerline{$^1$Department of Mathematics, Faculty of Science, Srinakharinwirot University, Bangkok 10110, Thailand
% address goes in here
}
\centerline{\texttt{$^1$thitarie@swu.ac.th}
% email goes in here
}

% email goes in here
%\support{Include acknowledgement of support here}

\bigskip

\begin{abstract}
Finding necessary and sufficient conditions for isomorphism between two semigroups of order-preserving transformations   over an infinite domain with restricted range was an open problem in \cite{FHQS}.  In this paper, we show a proof strategy to answer that question. 
\end{abstract}

\bigskip
\section{\bf Introduction}

For a nonempty set $X$, let $T(X)$ be the full transformation semigroup under composition of all maps from $X$ to $X$. When $X$ is a partially ordered set (poset), a mapping  $\a$ in $T(X)$ is called \emph{order-preserving}\; if  $x\leq y$ implies $x\a\leq y\a$  for all $x,y\in X$,\; and  $\a$ is \emph{regressive}\; if $x\a\leq x$ for all $x\in X$. We denote by  $T_{OP}(X)$ and $T_{RE}(X)$ the subsemigroups of $T(X)$ of  all order-preserving maps and all regressive maps on $X$, respectively. The semigroups of order-preserving maps  was first introduced by Howie in \cite{Ho}.   \\
\indent For $\a \in T(X)$, let $\ran \a$ denote the range of $\a$.  In 1975, Symons \cite{S} introduced and studied the subsemigroup $T(X,X')$  where $\emptyset\neq X' \subseteq X $  of $T(X)$ consisting of $\a\in T(X)$ with  $\ran \a \subseteq X'$.  Subsemigroups of transformations (with restricted range) of $T(X)$ of this type have been studied extensively, including our work which we will mention later on. Regarding the semigroups of regressive type,  in 1996,  Umar  proved that for any chains $X$ and $Y$, $T_{RE}(X)\cong T_{RE}(Y)$ if and only if $X$ and $Y$ are order-isomorphic (see in \cite{U2}).  Later in \cite{SAK}, T. Saito, et al.  generalized this result to partially ordered sets. They introduced  the adjusted partially ordered set $A(X)$ of a poset $X$ and proved that the order-isomorphism between $A(X)$ and $A(Y)$ is a necessary and sufficient condition  for the two semigroups to be isomorphic.\\
 \indent In this paper, we are also interested in studying the isomorphisms of  subsemigroups of transformations with restricted range. Now, let us  introduce the subsemigroups which will be of particular interest to us in this paper.\\
\indent For  a partially ordered set $X$ and a subset  $X'$ of $X$, we let 
\begin{align*}
T_{OP}(X,X')&:=T_{OP}(X)\cap T(X,X'),\\
T_{RE}(X,X')&:=T_{RE}(X)\cap T(X,X').
\end{align*}
Then both of these are subsemigroups of $T(X,X')$. \\
\indent In 2012, Udomkavanich and Jitjankarn  proved in  \cite{UJ2} that $T_{RE}(X,X')\cong T_{RE}(Y,Y')$ if and only if two adjusted chains $\mathcal{A}(X,X')$ and $\mathcal{A}(Y,Y')$ are order-structural isomorphic. This result leads us to study the isomorphism theorems for the semigroups of order-preserving type. It is known (e.g.,\cite{Ly}, page 222-223) that for  posets $X$ and $Y$, $T_{OP}(X)\cong T_{OP}(Y)$ if and only if $X$ and $Y$ are either order-isomorphic or order-anti-isomorphic. These necessary and sufficient conditions also hold  for the isomorphisms on the semigroups of  partial order-preserving  transformations (see in \cite{KMR}). In 2013, Fernandes, et al. \cite{FHQS} show that these   conditions apply for  $T_{OP}(X,X')$ and $T_{OP}(Y,Y')$ to be  isomorphic when $X$ and $Y$ are finite  as well. In this paper, we study the case when $X$ and $Y$ are infinite chains. Since $T_{OP}(X,X')$ is trivial when $|X'| = 1$, we omit this case. 

\indent Throughout the paper, we assume that $X$ and $Y$ are chains,  $|X'| > 1$, and $ |Y'| >1$. The following statement is known.
\begin{align}
\text{If there is an order-(anti)-isomorphism} \;\t:X\rightarrow Y\;\;\text{such that}\;\notag\\
(X')\t=Y'\;\;\;\text{for some}\;X'\subseteq X\;\text{and}\;Y'\subseteq Y,\quad\quad\quad \label{eq101}\\
\text{then}\;\; T_{OP}(X,X')\cong T_{OP}(Y,Y').\quad\quad\quad\quad\quad\quad\quad\quad\quad\quad\quad\quad\quad\notag
\end{align}
\indent  It is natural to ask whether the converse of the above result holds. Nevertheless, our work shows that it may not be the case if $|X'| = 2$. To be precise, we derive that the converse of the statement (\ref{eq101}) holds  when $|X'|\geq3$. 

\noindent To prove the statements, we apply in a similar fashion to  \cite{UJ2} the idea of  using adjusted chains. To do so, we will first introduce some notation and definitions that will be useful in Section 2. In Section 3, some homomorphism properties which are preserved under  isomorphism will be given. Lastly, the  isomorphism theorems for the semigroups of the type $T_{OP}(X,X')$ when $X$  is an infinite chain are determined in Section 4. \\

%%%%%%%%%%%%%%%%%%%%%%%%%%%%%%%%
\section{\bf Basic notations and results}
%%%%%%%%%%%%%%%%%%%%%%%%%%%%%%%%. 
\indent  Let  $\C'$ be a subchain of a  chain $\C$. Let $\{\C\backslash{\C'}\}$ denote the  set of all equivalence classes of  $\C\backslash{\C'}$ such that each class contains all elements in $\C\backslash{\C'}$ with no elements in $\C'$ lying between them. Then we consider $\{\C\backslash{\C'}\}\cup\C'$ as a chain under the partial order induced by the chain $\C$ in the natural way. This chain is an adjusted chain, denoted by $\mathcal{A}\{\C,\C'\}$.\\
\indent For each $a,b\in\C$ with $a<b$, the intervals $(a,b), \;[a,b), \;(a,b], \;[a,b]$ in $\C$ are defined naturally and we define the following intervals.
 $$\begin{array}{rclrcl}
(\leftarrow a]&:=&\{z\in \C~:~z\leq a\},& [a\rightarrow )&:=&\{z\in \C~:~z\geq a\},\\
(\leftarrow a)&:=&\{z\in \C~:~z<a\},& (a\rightarrow )&:=&\{z\in \C~:~z> a\}.
\end{array}$$
 
\indent For a nonempty subset $V$ of a chain $\C$, $V$ is said to be \emph{convex} if for $x, y, z \in \C$ such that $x\leq z\leq y$,  $x,y \in V$ implies  $z\in V$;  $V$  is called an \emph{upper$($lower$)$-convex subset} of $\C$ if\; $x<y$ ($x>y$)\; for all $x\in\C\backslash{V}$ and $y\in V$.\\ 
\indent For a convex subset $V$ of $\C$, we define  
$$\begin{array}{rclrcl}
%(\leftarrow [c]]&=&\{z\in \C~:~z\leq [c]\},& [[c]\rightarrow )&=&\{z\in \C~:~z\geq [c]\},\\
(\leftarrow V)&:=&\{z\in \C\backslash V~:~z \; \textnormal{is a lower bound of $V$ in $\C$}\},&\\ (V\rightarrow )&:=&\{z\in \C\backslash V~:~z\; \textnormal{is an upper bound of $V$ in $\C$}\}.&
\end{array}$$
 
For convenience, if $a\in\C$, let $\C_a$ be the element of $T(\C)$ whose range is $\{a\}$.\\
\indent Given $[k]\in\{\C\backslash{\C'\}}$, $k\in \C\backslash{\C'}$. We will define some order-preserving maps of  $T(\C,\C')$ as follows:% For $a,b,c,d\in X'$,
\begin{itemize}
\item For a convex subset $A$ of $[k]$ and $a,b,c\in \C'$ such that $a\leq b<[k]<c$ (or $a<[k]<b\leq c$), we write\smallskip\\
 \indent\quad\quad\quad\quad\quad\quad \quad\quad\quad\quad
  $\omega_{_{a:A_b:c}}:={\scriptsize\begin{pmatrix}
(\leftarrow A)& A& (A\rightarrow)\\
a&b&c
\end{pmatrix}}$  
\smallskip\smallskip\\
where $\ran(\omega_{_{a:A_b:c}}) =\{a,b,c\} $ and $b(\omega_{_{a:A_b:c}})^{-1}=A$ if $a\neq b$. 

\item When $[k]=\min\mathcal{A}\{X\backslash{X'\}}$,\; for a lower-convex subset $L$ of $[k]$ \;and $a,b\in \C'$ such that  \;$[k]<a<b$, we write\smallskip\\
 \indent\quad\quad\quad\quad\quad\quad \quad\quad\quad\quad\quad
  $\omega_{_{L_a:b}}:={\scriptsize\begin{pmatrix}
L&  (L\rightarrow)\\
a&b
\end{pmatrix}}$
\smallskip\smallskip\\
where $\ran(\omega_{_{L_a:b}}) =\{a,b\} $ and $a(\omega_{_{L_a:b}})^{-1}=L$. 

\item When $[k]=\max\mathcal{A}\{X\backslash{X'\}}$,\;for an upper-convex subset $U$ of $[k]$ \;and $a,b\in \C'$ such that  \;$a<b<[k]$, we write\smallskip\\
 \indent\quad\quad\quad\quad\quad\quad \quad\quad\quad\quad\quad
  $\omega_{_{a:U_b}}:={\scriptsize\begin{pmatrix}
(\leftarrow U)&  U\\
a&b
\end{pmatrix}}$
\smallskip\smallskip\\
where $\ran(\omega_{_{a:U_b}}) =\{a,b\} $ and $b(\omega_{_{a:U_b}})^{-1}=U$. 
\end{itemize}

For $\a\in T(\C)$, we denote $Fix(\a)=\{x\in \C ~:~x\a=x\}$. \\
\indent
For $\a\in T(\C,\C')$, we define \emph{the partial graph of transformation} $\a$, denoted by $\Gamma_\a:=(\C',\ran\a,E_\a)$, in the following way: $\C'$ is the set of upper vertices,  $\ran\a$ is the set of lower vertices such that all vertices are placed in order, and $E_\a$ is the set of (directed) edges which each element is in the form  $x\a y$, where $x\a=y$ for $x,y \in \C'$. Notice that the number of  components in each partial graph is equal to the number of elements in its range. Furthermore,   the   components,   considered from  left to right, are placed in the  same order as their related elements in the range.

\noindent\emph{Example.} For the transformation $\a\in  T(\{1,2,\ldots,9\},\{1,3,5,7,9\})$ defined by
$$\a={\scriptsize\begin{pmatrix}
1& 2 &3 &4&5&6&7&8&9\\
1& 5 &1 &9&5&5&5&3&5
\end{pmatrix}},$$
the set of  upper vertices is $\{1,3, 5, 7, 9\}$, the set of lower vertices is $\{1,3,5,9\}$ and $E_\a = \{1\a 1, 3\a 1, 5\a 5, 7\a 5, 9\a 5\}$. Then the graph $\Gamma_\a$ has the following form:

$$
{\scriptsize\xygraph{
!{<0cm,0cm>;<0.6cm,0cm>:<0cm,0.6cm>::}
!{(.5,.5) }*{\bullet}="y"     !{(1.5,.5) }*{\bullet}="z"      !{(2.5,.5)}*{\bullet}="d"    !{(3.5,.5)}*{\bullet}="g"    !{(4.5,.5)}*{\bullet}="i" %%%%%%%%%%%%%%
!{(.5,-.5) }*{\bullet}="w"    !{(1.5,-.5) }*{\bullet}="aa"                                                                                                                                        !{(2.5,-.5)}*{\bullet}="f"   !{(4.5,-.5)}*{\bullet}="fff"    
"y"-"w"  "z"-"w"  "d"-"f"  "g"-"f"  "i"-"f" 
}} \ \ $$ 

\bigskip
\noindent The partial graph $\Gamma_\a$ has four components placed in order from left to right.

%%%%%%%%%%%%            Thm  2.1
\begin{thm}\label{thm2.1} If $T_{OP}(X,X')\cong T_{OP}(Y,Y')$, then   $ X'$ and $Y'$ are either order-isomorphic or order-anti-isomorphic. 
\end{thm}
\begin{proof}
Let $\v:T_{OP}(X,X')\rightarrow  T_{OP}(Y,Y')$ be an isomorphism. For each $a\in X'$, there is an element $\bar{a}\in Y'$ such that  $(X_a)\v=Y_{\bar{a}}$ by idempotent and right zero properties of $X_a$ and $Y_{\bar{a}}$. The map $a\mapsto\bar{a}$ becomes a bijective map from $X'$ onto $Y'$. It remains to show that this map is either order-preserving or order-anti-preserving. Let $a,b,s,t\in X'$ be such  that $a<b$ and $s<t$. Since $Y'$ is a chain and the map is one-to-one, it must be that  $\bar{a} < \bar{b}$  or  $\bar{a} > \bar{b}$, and  $\bar{s} < \bar{t}$  or  $\bar{s} > \bar{t}$. Now, we have ${\scriptsize
\begin{pmatrix}
(\leftarrow a]& (a\rightarrow)\\
s&t \end{pmatrix}}\in T_{OP}(X,X')$ such that 
 $$X_a{\scriptsize
\begin{pmatrix}
(\leftarrow a]& (a\rightarrow)\\
s&t \end{pmatrix}}=X_s~\;\text{and}~\; X_b{\scriptsize
\begin{pmatrix}
(\leftarrow a]& (a\rightarrow)\\
s&t \end{pmatrix}}=X_t.   $$Then
$$Y_{\bar{a}}\left({\scriptsize
\begin{pmatrix}
(\leftarrow a]& (a\rightarrow)\\
s&t \end{pmatrix}}\v\right)=Y_{\bar{s}}~\;\text{and}\;~ Y_{\bar{b}}\left({\scriptsize
\begin{pmatrix}
(\leftarrow a]& (a\rightarrow)\\
s&t \end{pmatrix}}\v\right)=Y_{\bar{t}}.$$
Consequently, \begin{center}$\bar{a}{\scriptsize
\begin{pmatrix}
(\leftarrow a]& (a\rightarrow)\\
s&t \end{pmatrix}}\v=\bar{s}$ \;and\; $\bar{b}{\scriptsize
\begin{pmatrix}
(\leftarrow a]& (a\rightarrow)\\
s&t \end{pmatrix}}\v=\bar{t}$.\end{center} 

Since ${\scriptsize
\begin{pmatrix}
(\leftarrow a]& (a\rightarrow)\\
s&t \end{pmatrix}}\v\in T_{OP}(Y,Y')$, it follows that  $\bar{a}<\bar{b}$ implies  $\bar{s}<\bar{t}$   and  $\bar{a}>\bar{b}$ implies $\bar{s}>\bar{t}$. This proves that $ X'$ and $Y'$ are either order-isomorphic or order-anti-isomorphic.
\end{proof}
\indent From now on,  let   $\varphi$ denote an isomorphism from  $T_{OP}(X,X')$ and $T_{OP}(Y,Y')$. The order-(anti)-isomorphism from $X'$ onto $Y'$, defined in the proof of Theorem \ref{thm2.1},  is denoted by $\t_\v$. It is easy to see that  the order-(anti)-isomorphism  $\t_{\v^{-1}}$ from $Y'$ onto $X'$, induced by the isomorphism $\v^{-1}$, is the inverse function of $\t_\v$. That is, 
\[\t_{\v^{-1}}=(\t_\v)^{-1}.\]
Notice that by considering $\v^{-1}$ and $\t_{\v^{-1}}$ instead of $\v$ and $\t_\v$, respectively,  all results that hold for $\v$ also hold for $\v^{-1}$. 

%%%%%%%%%%%%%%%%%%%%%%%%%%%%%%%%
\section{\bf Some homomorphism properties}
%%%%%%%%%%%%%%%%%%%%%%%%%%%%%%%%. 
In this section, we study some properties of transformations which will be preserved under a homomorphism.   First, we will study the structure of $\a$ and $\a\v$ through $\t_\v$ when $\a\in T_{OP}(X,X')$. Then we derive that   two graphs of $\Gamma_\a$ and $\Gamma_{(\a)\varphi}$  are isomorphic. Moreover, the order of components (in the sense of partial graph)   is also preserved.  

\indent Without loss of generality,   we assume that $\t_\v$ is  order-preserving from now on. The other case that $\t_\v$ is order-anti-preserving can be done by the same  process. 

 %%%%%%%%%%%%%%%%%%%%%%%%%%%%%%%%%%%%%%%%%%%%%%%%%%%%%%%%%%%%%%0000000000000000000%%%%%%%%%%%%%%%%%%%%%%%
\begin{lemma}\label{lemma3.1} For each $\a\in T_{OP}(X,X')$, the following statements hold:
\begin{enumerate}[(i)]
\item  $(Fix(\a))\t_\v=Fix(\a\v)$.
\item  For $a\in\ran\a$ such that $a\a^{-1}\cap X'\neq\emptyset$,\\
\indent \ \ \ \ \ \ $\bar{a}\in\ran(\a\v)$ and $\bar{a}(\a\v)^{-1}\cap Y'=(a\a^{-1}\cap X')\t_\v$.
\end{enumerate}
In particular,   if $\a$ is an idempotent, then $(\ran\a)\t_\v=\ran(\a\v)$.      
\end{lemma}
%ooooooooooooooooooooooooooooooooooooooooooooooooooooooooooooooooooooooooooooooooooooooooooooooooooooooooooooo0000000000000000oooooooooooo
\begin{proof} (i) Let $a\in Fix(\a)$. Then $a\a=a$. Since $Y_{\bar{a}}(\a\v)=(X_a\v)(\a\v)=(X_a\a)\v=X_a\v=Y_{\bar{a}}$, it follows that   $ \bar{a}(\a\v)=\bar{a}=a\t_\v\in Fix(\a\v)$. Similarly, if   $\bar{s}\in Fix(\a\v)$, then $X_s\a=(Y_{\bar{s}}\v^{-1})\a=(Y_{\bar{s}}(\a\v))\v^{-1}=(Y_{\bar{s}})\v^{-1}=X_s$, that is, $s\a=s$. Then $\bar{s}=s\t_\v\in(Fix(\a))\t_\v$. \\
\indent(ii) For  $a\in\ran\a$  such that $a\a^{-1}\cap X'\neq\emptyset$, let $x\in a\a^{-1}\cap X'$.   Then 
$a\in Fix\left(X_x\a\right)$, by (i),  $\bar{a}\in Fix\left((X_x\v)(\a\v)\right)$. That is, $\bar{a}\in\ran(\a\v)$. Since $\bar{x}(\a\v)=\bar{a}(X_x\v)(\a\v)=\bar{a}$, it follows that $\bar{x}\in \bar{a}(\a\v)^{-1}\cap Y'$. Then $(a\a^{-1}\cap X')\t_\v\subseteq\bar{a}(\a\v)^{-1}\cap Y'$. Similarly, by considering $\v^{-1}$ instead of $\v$,  $\bar{a}(\a\v)^{-1}\cap Y'\subseteq(a\a^{-1}\cap X')\t_\v$. Thus the  equality is obtained. 
\end{proof}
%%%%%%%%%%%%%%%%%%%mmmmmmmmmmmmmmmmmmmmmmmmmmmmmmmmmmmmmmmmmmmmmmmmmmmmm%%2.3
\begin{lemma}\label{lemma3.2}For each $\a\in T_{OP}(X,X')$,
 if $b\in\ran\a$ and $b\a^{-1}\cap X'=\emptyset$, then $\bar{b}\in\ran(\a\v)$.
\end{lemma}
%%%%%%%%%%%%%%%%%%%mmmmmmmmmmmmmmmmmmmmmmmmmmmmmmmmmooooooooooooooooooooooooooooooooooooooooooooooooooooooooooooooooooooooooooooooooooommmmmmmmmmmmm2.3
\begin{proof}  Let $b\in\ran\a$ and $b\a^{-1}\cap X'=\emptyset$. Assume that $b$ is neither maximum nor minimum in $X'$. Choose $a,c\in X'$ such that $a<b<c$ and let  $\epsilon_b={\scriptsize\begin{pmatrix}
(\leftarrow b)&b& (b\rightarrow) \\
a& b& c
\end{pmatrix}}$. 
 Then $\epsilon_b$ is an idempotent with $b(\epsilon_b)^{-1}\cap X' =\{b\}$. By Lemma \ref{lemma3.1}, $\bar{b}(\epsilon_b\v)^{-1}\cap Y'=\{\bar{b}\}$. Suppose in the contrary that $\bar{b}\notin \ran(\a\v)$. Then we have $\bar{b}\notin\ran ((\a\v)(\epsilon_b\v))$.  Since $|\ran(\a\v)(\epsilon_b\v)|$ is finite, this guarantees the existence of an idempotent $\mu$ in $T_{OP}(Y,Y')$ with $\ran\mu=\ran(\a\v)(\epsilon_b\v)$. Then $\mu\v^{-1}$ is an idempotent in $T_{OP}(X,X')$, by Lemma \ref{lemma3.1}, $b\notin \ran(\mu\v^{-1})$. However,
 $$\a\epsilon_b(\mu\v^{-1})=((\a\v)(\epsilon_b\v)\mu)\v^{-1}=((\a\v)(\epsilon_b\v))\v^{-1}=\a\epsilon_b,$$   
which is a contradiction. If $b$ is either  maximum or minimum, it can be proved in the same way   by defining $\epsilon_b$ as before and choosing $a=b$ if $b$ is minimum, and  $c=b$ if $b$ is maximum. 
\end{proof}
%%%%%%%%%%%%%%%%%%%%%%%%%Proposition  2.4 %%%%%%%%%%%%%**************************************************%%%%%%%%%
By  Lemma \ref{lemma3.1} and \ref{lemma3.2},  the following proposition is directly obtained.

\begin{prop}\label{prop3.3}For each $\a\in T_{OP}(X,X')$, we have 
\begin{enumerate}[(i)]
\item  $(\ran\a)\t_\v=\ran(\a\v)$. 
\item  For any $a\in\ran\a$, $\bar{a}(\a\v)^{-1}\cap Y'=(a\a^{-1}\cap X')\t_\v$.
\end{enumerate}
\end{prop}

\noindent This proposition leads us to define an interesting equivalence relation on the semigroup of full transformations with restricted range. 

Given a transformation $\a:X\rightarrow X'$, the $\a$-\emph{structure} is  the partial graph $\Gamma_\a$ and putting its components in the same order as their related elements in the range.
Here we define an equivalence relation $\mathcal{K}$ on  $T(X,X')$ by 
$$\a\mathcal{K}\b \ \ \ \  \text{iff} \ \  \ \ \a\text{-structure} \ \  \text{and} \ \ \b\text{-structure} \ \ \text{are the same}.$$
Indeed, it is equivalent to $\a|_{_{X'}}=\b|_{_{X'}}$ and $\ran{\a}=\ran{\b}$.
The $\K$-class containing $\a$ is denoted by $\K_\a$.
It is very clear that when $X'=X$, $T(X,X)$ is $\K$-trivial. By  Proposition \ref{prop3.3}, we have that  $\K_\a$ and $(\K_\a)\v=\K_{\a\v}$ have the same structure for all $\a\in T_{OP}(X,X')$.

\indent Next, we will construct an extension of $\theta_{\varphi}$ to be an order-isomorphism on the adjusted chains.

\begin{lemma}\label{lemma3.4}Suppose that two classes $[k_1]$ and $[k_2]$ are the minimum and the maximum  of $\mathcal{A}\{X,X'\}$, respectively.  Let $a,b\in X'$ be such that  $a<b$, and  $A\subseteq[k_1]$ and $B\subseteq[k_2]$ as a lower-convex subset and an upper-convex subset of $[k_1]$ and $[k_2]$, respectively.  Then
\[( \omega_{_{{A}_a:b}}) \v=\omega_{_{{C}_{\bar{a}}:\bar {b}}}\quad\quad\text{and}\quad\quad( \omega_{_{a:{B}_b}}) \v=\omega_{_{\bar {a}:{D}_{\bar{b}}}}\]
for some lower-convex  $C$ and upper-convex  $D$ of the minimum and  the maximum   of $\{Y\backslash{Y'\}}$, respectively. 
\end{lemma}

\begin{prop} \label{prop3.5} For each $[k]\in\{X\backslash{X'}\}$, there is a corresponding $[t_k]\in\{Y\backslash{Y'\}}$ such 
that the extended map of $\t_\v$ from $X'\cup\{[k]\}$ onto $Y'\cup\{[t_k]\}$ is an order-isomorphism. Moreover, $|[k]|=|[t_k]|$.
\end{prop}
\begin{proof}Let $[k]\in \{X\backslash{X'}\}$ be such that 
$a<[k]<b$  for some $a,b\in X'$. We choose $\omega_{_{a:A_a:b}}$ as an  idempotent in $T_{OP}(X,X')$ whose  range is $\{a,b\}$. Since two partial graphs of transformations  
$\Gamma_{ \omega_{_{a:A_a:b}}}$ and 
$\Gamma_{(\omega_{_{a:A_a:b}})\v}$ have the same structure, by 
 Proposition \ref{prop3.3}, it follows that 
 $(\K_{\omega_{_{a:A_a:b}}})\v=\K_{(\omega_{_{a:A_a:b}})\v}$. Due to the  structure of $\Gamma_{ \omega_{_{a:A_a:b}}}$, the cardinality of $\K_{\omega_{_{a:A_a:b}}}$ is depending only on $[k]$. Indeed, $|\K_{\omega_{_{a:A_a:b}}}|=|[k]|$ . This imply the  existence of   
$[t_k]\in\{Y\backslash{Y'}\}$ with $\bar{a}<[t_k]<\bar{b}$) and  
$|\K_{(\omega_{_{a:A_a:b}})\v}|=|[t_k]|$.\\
\indent Suppose  $[k]$ is  maximum $($or minimum$)$ in $\mathcal{A}\{X,X'\}$.  For any $a,b\in X'$ such that $a<b$, we consider $\omega_{_{a:[k]_b}}$ $($or  $\omega_{_{[k]_a:b}}$ $)$. By Lemma \ref{lemma3.4} and using the same argument, our proof is finished. 
\end{proof}

From Proposition \ref{prop3.5}, the union of all these extensions form  an order-isomorphism, denoted by $\widehat{\t_\v}$ (with respect to $\t_\v$), from $\mathcal{A}\{X,X'\}$ onto $\mathcal{A}\{Y,Y'\}$ such that
$$x\mapsto \bar{x}\;\;\text{and}\;\;[k]\mapsto[t_k]$$
for $x\in X'$ and $[k]\in\{X\backslash{X'\}}$. We notice that $\widehat{\t_\v}$ is an order-structural isomorphism (as defined in \cite{UJ2}).   This conclusion results in the isomorphism theorems between the two semigroups for an infinite discrete chain.

\begin{thm}\label{thm3.6} Let $X$ and $Y$ be  infinite discrete chains. Then  $T_{OP}(X,X')\cong T_{OP}(Y,Y')$ if and only if there is an order-(anti)-isomorphism $\t:X\rightarrow Y$ such that $(X')\t=Y'$.
\end{thm}

 Nevertheless, the property  that  $|[k]|=|[t_k]|=|([k])\widehat{\t_\v}|$ is not sufficient to determine  the isomorphism for an uncountable chain. As a result, we study more of homomorphism properties associated with a class of $\{X\backslash{X'}\}$.

\begin{lemma}\label{lemma3.6}Let  $[k]\in \{X\backslash{X'\}}$ be such that $a<b<[k]<c$ $($or $a<[k]<b<c$$)$ for some $a,b,c\in X'$. Then for each convex subset $A$ of $[k]$, 
\[(\omega_{_{a:A_b:c}}) \v=\omega_{_{\bar{a}:{B}_{\bar{b}}:\bar {c}}}\]
for some convex  subset $B$ of  $[t_k]\in \{Y\backslash{Y'\}}$ with  $\bar{a}<\bar{b}<[t_k]<\bar{c}$ $($or $\bar{a}<[t_k]<\bar{b}<\bar{c}$$)$.
\end{lemma}
\begin{proof} By Proposition \ref{prop3.3}, it follows that  $\ran((\omega_{_{a:A_b:c}})\v)=\{\bar{a},\bar{b},\bar{c}\}$. Since $(a\omega_{_{a:A_b:c}}^{-1}\cup \;c\omega_{_{a:A_b:c}}^{-1})\cap X'=X'$, by Lemma \ref{lemma3.1}, we have that $(\bar{a}((\omega_{_{a:A_b:c}})\v)^{-1}\cup \bar{c}((\omega_{_{a:A_b:c}})\v)^{-1})\cap Y'=Y'$. As $((\omega_{_{a:A_b:c}})\v)$ is  order-preserving such that $\bar{b}$ is in its range, there exists the unique class in $\{Y\backslash{Y'\}}$, namely $[t_k]$, containing a convex subset $\bar{b}((\omega_{_{a:A_b:c}})\v)^{-1}$.
\end{proof}
%%%%&&&&&&&&66666666666
\begin{prop}\label{prop3.8}Let $[k]\in \{X\backslash{X'\}}$ be such that $a<b<[k]<c$ $($or $a<[k]<b<c$$)$ for some $a,b,c\in X'$. Then for each $x\in[k]$, 
\[(\omega_{_{a:\{x\}_b:c}}) \v=\omega_{_{\bar{a}:\{y\}_{\bar{b}}:\bar {c}}}\]
for some $y\in[t_k]$.
\end{prop}
\begin{proof} Let $f$ and $g$ stand for two idempotents in $T_{OP}(X,X')$ such that $\ran f=\ran g=\{a,c\}$ with $\{b,c\}\subseteq cf^{-1}$  and $\{a,b\}\subseteq ag^{-1}$. Let $B$ be a convex subset of $[k]$ such that $|B|>1$. By  Lemma \ref{lemma3.6}, we obtain that $(\omega_{_{a:B_b:c}}) \v=\omega_{_{\bar{a}:M_{\bar{b}}:\bar {c}}}$ for some convex subset $M$ of $[t_k]$.  Suppose  in the contrary that $M=\{y\}$.   We choose $L$ and $U$ as two convex subsets of $B$ which form a partition of $B,$ and $L$ is a lower bound of $U$.  Since $\omega_{_{a:{U}_{b}:c}}g=\omega_{_{a:{B}_{b}:c}}g$, it follows that
\begin{equation*}
(\omega_{_{a:{U}_{b}:c}}\v)(g\v)=(\omega_{_{a:{B}_{b}:c}}\v)(g\v)=\omega_{_{\bar{a}:\{y\}_{\bar{b}}:\bar{c}}}(g\v).
\end{equation*}
Then $\bar{b}(\omega_{_{a:{U}_{b}:c}}\v)^{-1}$ is an upper-convex subset of $(\leftarrow y]$. Since $\omega_{_{a:{L}_{b}:c}}f=\omega_{_{a:{B}_{b}:c}}f$, we have
\begin{equation*}
(\omega_{_{a:{L}_{b}:c}}\v)(f\v)=(\omega_{_{a:{B}_{b}:c}}\v)(f\v)=\omega_{_{\bar{a}:\{y\}_{\bar{b}}:\bar{c}}}(f\v).
\end{equation*}
Then  $\bar{b}(\omega_{_{a:{L}_{b}:c}}\v)^{-1}$ is a lower-convex subset of $[y\rightarrow )$. It can be seen that $\omega_{_{a:{L}_{b}:c}}g=\omega_{_{a:{U}_{b}:c}}f$. Then $(\omega_{_{a:{L}_{b}:c}}\v)(g\v)=(\omega_{_{a:{U}_{b}:c}}\v)(f\v)$ which contradicts to
\begin{equation*}
\bar{a}=\bar{b}(g\v)=y(\omega_{_{a:{L}_{b}:c}}\v)(g\v)=y(\omega_{_{a:{U}_{b}:c}}\v)(f\v)=\bar{b}(f\v)=\bar{c}.
\end{equation*}
\end{proof}

%%%%%%%%+++++++++++++++++++++++++ Prop3.9
\begin{prop}\label{prop3.9} For  $a,b,c\in X'$ with $a<b<c$, the following statements hold:
\begin{enumerate}[(i)]
\item  If  $[k]=\max\mathcal{A}\{X\backslash{X'}\}$, then for $x\in[k]$, \\\indent \quad\quad \quad\quad  \quad\quad  \quad\quad 
          \;$(\omega_{_{a:[x\rightarrow)_{_c}}})\v=\omega_{_{\bar{a}:[y\rightarrow)_{_{\bar{c}}}}}$\quad \quad for some $y\in [t_k]$.
\item  If  $[k]=\min\mathcal{A}\{X\backslash{X'}\}$, then for $x\in[k]$,  \\\indent \quad\quad \quad\quad  \quad\quad  \quad\quad
          \;$(\omega_{_{(\leftarrow x]_{_a}:c}})\v=\omega_{_{(\leftarrow y]_{_{\bar{a}}}:\bar{c}}}$\quad \quad for some $y\in [t_k]$.
\end{enumerate}
\end{prop}
\begin{proof}(i) Suppose  $[k]=\max\mathcal{A}\{X\backslash{X'}\}$. Let $x\in[k]$. Suppose that $(x\rightarrow)\neq\emptyset$. We let $\a={\scriptsize 
\begin{pmatrix}
(\leftarrow x)& x &(x\rightarrow) \\
a&b & c 
\end{pmatrix}}$ and  $\b={\scriptsize 
\begin{pmatrix}
(\leftarrow b)& [b\rightarrow) \\
a&  c 
\end{pmatrix}}$. Clearly, 
$\a\b=\omega_{_{a:[x\rightarrow)_{_c}}}$.
Then $(\a\v)(\b\v)=(\omega_{_{a:[x\rightarrow)_{_c}}})\v$. By applying the same process as in the proof of Proposition \ref{prop3.8}, we obtain that $|\bar{b}(\a\v)^{-1}|=1$. Since
$\bar{c}(\omega_{_{a:[x\rightarrow)_{_c}}}\v)^{-1}=\bar{b}(\a\v)^{-1}\dot{\cup}\bar{c}(\b\v)^{-1}$
where $\bar{b}(\a\v)^{-1}$ is a lower-convex subset of $\bar{c}(\omega_{_{a:[x\rightarrow)_{_c}}}\v)^{-1}$.
These imply that $\bar{c}(\omega_{_{a:[x\rightarrow)_{_c}}}\v)^{-1}=[y\rightarrow)$ for some $y\in[t_k]$.\\
\indent (ii) can be proved similarly to (i).
\end{proof}

%%%%%%%%%%%%%%%%%%%3333333333333333333333333333333333333
\section{\bf Isomorphism theorems}
%Recall that $\t_v$ is order-isomorphic.\\

In the last section, we  take care of the case $|X'|=2$. 
For convenience,  we here denote $T_{OP}(X,X')$ by $\mathcal{O}[_{M_1}1_{M_2}2_{M_3}]$ where $M_1, M_2$ and $M_3$ are three classes in $\{X\backslash{X'}\}$ .  We observe that there are only 5 classes in $\mathcal{O}[_{M_1}1_{M_2}2_{M_3}]/\mathcal{K}$ whose  partial graph of transformations is one of the  following forms:\begin{center}
 \bigskip
  $\Gamma_{\lambda_1} :\;
{\scriptsize\xygraph{
!{<0cm,0cm>;<0.6cm,0cm>:<0cm,0.6cm>::}
!{(.5,.5) }*{\bullet}="y"          !{(2.5,.5)}*{\bullet}="d"       %%%%%%%%%%%%%%
!{(.5,-.5) }*{\bullet}="yy"             
"y"-"yy"  "d"-"yy"    
}}$  \quad \quad  
$\Gamma_{\lambda_2} :\;
{\scriptsize\xygraph{
!{<0cm,0cm>;<0.6cm,0cm>:<0cm,0.6cm>::}
!{(.5,.5) }*{\bullet}="y"          !{(2.5,.5)}*{\bullet}="d"         %%%%%%%%%%%%%%
                                   !{(2.5,-.5)}*{\bullet}="dd"         
"y"-"dd"  "d"-"dd"   
}}$   \quad \quad  
$\Gamma_{\lambda_3} :\;
{\scriptsize\xygraph{
!{<0cm,0cm>;<0.6cm,0cm>:<0cm,0.6cm>::}
!{(.5,.5) }*{\bullet}="y"          !{(2.5,.5)}*{\bullet}="d"         %%%%%%%%%%%%%%
!{(.5,-.5) }*{\bullet}="yy"        !{(2.5,-.5)}*{\bullet}="dd"          
"y"-"yy"  "d"-"dd" 
}}$\bigskip\\
 $\Gamma_{\lambda_4} :\;
{\scriptsize\xygraph{
!{<0cm,0cm>;<0.6cm,0cm>:<0cm,0.6cm>::}
!{(.5,.5) }*{\bullet}="y"          !{(2.5,.5)}*{\bullet}="d"         %%%%%%%%%%%%%%
!{(.5,-.5) }*{\bullet}="yy"         !{(2.5,-.5)}*{\bullet}="dd"          
"y"-"yy"  "d"-"yy" 
}}$ \quad \quad \quad
$\Gamma_{\lambda_5}:\;
{\scriptsize\xygraph{
!{<0cm,0cm>;<0.6cm,0cm>:<0cm,0.6cm>::}
!{(.5,.5) }*{\bullet}="y"          !{(2.5,.5)}*{\bullet}="d"         %%%%%%%%%%%%%%
!{(.5,-.5) }*{\bullet}="yy"        !{(2.5,-.5)}*{\bullet}="dd"         
"y"-"dd"  "d"-"dd"   
}}$  \bigskip
 
\end{center}
%%%%%%%%%%%%%%
The following results are directly derived.
%%%%%%%%%%%
\begin{lemma}\label{lemma4.1} For  $\mathcal{K}_{\lambda_i}\in\mathcal{O}[_{M_1}1_{M_2}2_{M_3}]/\mathcal{K}$, $(i=1,\ldots,5)$, we have that
\begin{enumerate}[(i)]
\item  $\mathcal{K}_{\lambda_1}$ and $\mathcal{K}_{\lambda_2}$ are trivial,
\item  
   $|\mathcal{K}_{\lambda_3} | = |M_2|+1$,
\item   \; $|\mathcal{K}_{\lambda_4}|=|M_3|$\; and\; $|\mathcal{K}_{\lambda_5}|=|M_1|$.
\end{enumerate}
\end{lemma}
%%%%%%%%%%
\begin{proof} Since there are only two constant maps, (i)  is proved. To show (ii), it is easy to see that each element in $M_2$ determine the consequent map in $\mathcal{K}_{\lambda_3}$ and vice versa. Hence the bijection between the two sets is constructed. The same idea can also be applied to show\; $|\mathcal{K}_{\lambda_4}|=|M_3|$\; and\; $|\mathcal{K}_{\lambda_5}|=|M_1|$.  
\end{proof}  
%%%%%%%%%%%%%%%%
\begin{thm}\label{thm4.2} $\mathcal{O}[_{M_1}1_{M_2}2_{M_3}]\cong \mathcal{O}[_{N_1}\bar{1}_{N_2}\bar{2}_{N_3}]$ if and only if $|M_i|=|N_i|$ for all $i=1,2,3$.
\end{thm}
%%%%%%%%%%
\begin{proof} Suppose that for $i=1,\ldots,5$,  $\lambda_i$ and $\gamma_i$ are two representations of order-preserving maps having the same partial graph in $\mathcal{O}[_{M_1}1_{M_2}2_{M_3}]$ and  $  \mathcal{O}[_{N_1}\bar{1}_{N_2}\bar{2}_{N_3}]$, respectively.  By Lemma \ref{lemma4.1},  we let  $f_i$ be a bijection from $\mathcal{K}_{\lambda_i}$ onto $\mathcal{K}_{\gamma_i}$ for $i=1,\ldots,5$.   To show that $\v:=f_1\cup f_2\cup\cdots \cup f_5:\mathcal{O}[_{M_1}1_{M_2}2_{M_3}]\rightarrow \mathcal{O}[_{N_1}\bar{1}_{N_2}\bar{2}_{N_3}]$ is an isomorphism,  we let $\a\in\mathcal{O}[_{M_1}1_{M_2}2_{M_3}]$. It is easy to see that the pairwise composition of  five graph structures  can be one of the following
maps: for  $\b\in\mathcal{O}[_{M_1}1_{M_2}2_{M_3}],$  either $ \a\b=\lambda_1$, $\a\b=\lambda_2$ or $\a\b=\a,$
\begin{center}
\label{aggiungi}%\centering %
  \resizebox{4cm}{!} {\begin{tabular}{|c| c   c  c c c|}
    \hline
   & $\lambda_1 $& $\lambda_2$ & $\lambda_3$ &$\lambda_4 $&$ \lambda_5$ \smallskip\\ \hline
  $\lambda_1$   & $\lambda_1$      & $\lambda_2$     & $\lambda_1$       & $\lambda_1$    & $\lambda_2$      \smallskip \\ 
 $\lambda_2$   & $\lambda_1$      & $\lambda_2$     & $\lambda_2$       & $\lambda_1$    & $\lambda_2$      \smallskip \\  
 $\lambda_3$   & $\lambda_1$      & $\lambda_2$     & $\lambda_3$       & $\lambda_1$    & $\lambda_2$      \smallskip \\ 
$\lambda_4$    & $\lambda_1$      & $\lambda_2$     & $\lambda_4$       & $\lambda_1$    & $\lambda_2$      \smallskip \\ 
$\lambda_5$    & $\lambda_1$      & $\lambda_2$     & $\lambda_5$       & $\lambda_1$    & $\lambda_2$      \smallskip \\  \hline
  \end{tabular}}\smallskip\smallskip\\
\end{center} 
Suppose $\a\b=\lambda_1$.  One of the following statements hold:
\begin{center}
(i) \; $\b=\lambda_1$,  \quad\quad\quad
(ii) \; $\b\in\mathcal{K}_{\lambda_4}$,\quad\quad\quad
(iii) \; $\b\in\mathcal{K}_{\lambda_3}$ and $\a=\lambda_1$.\end{center}
It is clear that $(\a\v)(\b\v)=\gamma_1=(\lambda_1)\v$. \\
\indent For the rest, it can be proved  directly.  
\end{proof}  

%%%%%%%%%%%%%%%%
\noindent\emph{Example.}  Let $X = \mathbb{R}, X' = \{1,2\}, Y= [2,5), Y'= \{3,4\}$. 
Theorem \ref{thm4.2} tells us that $\mathcal{O}[_{(-\infty,1)}1_{(1,2)}2_{(2,\infty)}]\cong \mathcal{O}[_{[2,3)} 3_{(3,4)} 4_{(4,5)}]$, yet, it is clear that  $\mathbb{R}$ and $[2,5)$ are not order or order-anti-isomorphic.\\

%%%%%%%%%%%%%%%

Next, we will prove that when $|X'|\geq 3$, the converse of (\ref{eq101}) holds.
%%%%%%%%%%%%%%%%theorem4.3
\begin{thm}\label{thm4.3} Suppose that $|X'|\geq3$. Then
$T_{OP}(X,X')\cong T_{OP}(Y,Y')$ if and only if there is an order-$($anti$)$-isomorphism $\t$ from $X$ onto $Y$ such that  $(X')\t=Y'.$
\end{thm}
\begin{proof}It remains to show that for each $[k]\in \{X\backslash{X'}\}$, $[k]$ and $[t_k]$ are order-isomorphic. Let $[k]$ be a class in $\{X\backslash{X'}\}$. We will consider in two cases:\smallskip\smallskip\\
\textbf{Case 1.} $[k]=\max\mathcal{A}\{X\backslash{X'\}}$ or $\min\mathcal{A}\{X\backslash{X'\}}$.\\ WLOG, we assume that $[k]=\max\mathcal{A}\{X\backslash{X'\}}$. We choose $a,b,c\in X'$ with $a<b<c$. For any  $x,x'\in[k]$ with $x<x'$. Consider  $\omega_{_{a:[x\rightarrow)_{_c}}}$ and $\omega_{_{a:[x'\rightarrow)_{_c}}}$. By Proposition \ref{prop3.9},  we have
$(\omega_{_{a:[x\rightarrow)_{_c}}})\v=\omega_{_{\bar{a}:[y\rightarrow)_{_{\bar{c}}}}}$ and $(\omega_{_{a:[x'\rightarrow)_{_c}}})\v=\omega_{_{\bar{a}:[y'\rightarrow)_{_{\bar{c}}}}}$. 
%%%%%%%%%%%%%%%%%%%%%%%%%%%%%%%%%%%%%%
Let $\gamma={\scriptsize\begin{pmatrix}
(\leftarrow x)& [x,x') &[x'\rightarrow)\\
a&b&c
\end{pmatrix}}$. Then 
\begin{center}$\gamma{\scriptsize\begin{pmatrix}
(\leftarrow b)& [b\rightarrow)\\
a&c
\end{pmatrix}}=\omega_{_{a:[x\rightarrow)_{_c}}}$\;\; and\;\;
$\gamma{\scriptsize\begin{pmatrix}
(\leftarrow b]& (b\rightarrow)\\
a&c
\end{pmatrix}}=\omega_{_{a:[x'\rightarrow)_{_c}}}$.
\end{center}
It follows that 
\begin{center}$(\gamma\v){\scriptsize\begin{pmatrix}
(\leftarrow b)& [b\rightarrow)\\
a&c
\end{pmatrix}}\v=(\omega_{_{a:[x\rightarrow)_{_c}}})\v=\omega_{_{\bar{a}:[y\rightarrow)_{_{\bar{c}}}}}$, \\
$(\gamma\v){\scriptsize\begin{pmatrix}
(\leftarrow b]& (b\rightarrow)\\
a&c
\end{pmatrix}}\v=(\omega_{_{a:[x'\rightarrow)_{_c}}})\v=\omega_{_{\bar{a}:[y'\rightarrow)_{_{\bar{c}}}}}$.
\end{center}
Since $[y\rightarrow)=\bar{c}(\omega_{_{\bar{a}:[y\rightarrow)_{_{\bar{c}}}}})^{-1}=\bar{b}(\gamma\v)^{-1}\dot{\cup}\bar{c}(\gamma\v)^{-1}$ and
$[y'\rightarrow)=\bar{c}(\omega_{_{\bar{a}:[y'\rightarrow)_{_{\bar{c}}}}})^{-1}=\bar{c}(\gamma\v)^{-1}$, these imply that $y<y'$.
\smallskip\smallskip\\
\textbf{Case 2.} $[k]$ is neither $\max\mathcal{A}\{X\backslash{X'\}}$ nor $\min\mathcal{A}\{X\backslash{X'\}}$.\\ Then there are $a,b,c\in X'$ such that $a<b<[k]<c$ or $a<[k]<b<c$. \\By using Proposition \ref{prop3.8} and following the same proof  as in Case 1, we derive the result.
\end{proof}


\begin{thebibliography}{99}  

\bibitem{FHQS}
Fernandes V. H. , Honyam P. , Quinteiro T.M., Singha B. On semigroups of endomorphisms of a chain with restricted range. Semigroup Forum 2014;\;89: 77--104.


\bibitem{Ho}
Howie J.M. Products of idempotents in certain semigroups of order-preserving transformations. Proc. Edinburgh Math. Soc. 1971;\;17: 223--236.

\bibitem{KMR}
Kemprasit Y. , Mora W., Rungratgasame T.  Isomorphism theorems for semigroups of order-preserving  partial transformations. Int. J. Algebra 2010;\;4: 799--808.


\bibitem{Ly}
Lyapin E.S. Semigroups. Amer. Math. Soc., Providence, R.I., 1974. 


\bibitem{SAK}
  Saito T., Aoki K., Kajitori K.  Remarks on isomorphisms of regressive transformation semigroups. Semigroup Forum 1996;\;53: 129--134.

\bibitem{S} 
Symons J. S. V.  Some results concerning a transformation semigroup. J. Austral. Math. Soc. 1975;\; 19: 413--425.  


\bibitem{UJ2}
Udomkavanich P., Jitjankarn P.   A classification of regressive transformation semigroups on chains. Semigroup Forum 2012;\;85:  559--570.

\bibitem{U2}
Umar A.  Semigroups of order-decreasing transformations: The isomorphism theorem. Semigroup Forum 1996;\;53: 220--224.

\end{thebibliography}
\end{document}